\documentclass[12pt]{amsart}
\usepackage[all]{xy}
\usepackage{a4wide}
\usepackage{amssymb}
\usepackage{amsthm}
\usepackage{amsmath}
\usepackage{amscd}
\usepackage{verbatim}
\usepackage[all]{xy}
\usepackage{framed,color}
\usepackage{ytableau}

\usepackage{graphicx}
\usepackage{subfig}

\usepackage{amsmath,amsthm,amssymb}
\usepackage{enumitem}

\usepackage{float}
\usepackage{ifthen}
\usepackage{xargs}
\usepackage{url}
\usepackage[all]{xy}
\usepackage{xcolor}
\usepackage{comment}
\usepackage[obeyDraft,textsize=tiny]{todonotes}
\specialcomment{draftnote}{\vspace{1em}\begingroup\ttfamily\footnotesize\color{blue}\begin{minipage}{15cm}}{\end{minipage}\endgroup\vspace{1em}}
\specialcomment{alert}{\vspace{1em}\begingroup\ttfamily\footnotesize\color{blue}\begin{minipage}{15cm}}{\end{minipage}\endgroup\vspace{1em}}

\usepackage[refpage,intoc,noprefix]{nomencl}

\newcommand{\R}{\mathbb{R}} 
\newcommand{\C}{\mathbb{C}} 
\newcommand{\Z}{\mathbb{Z}} 
\newcommand{\N}{\mathbb{N}} 

\def\({\left(}
\def\){\right)}

\def\lp{\left(}
\def\rp{\right)}

\newtheorem{theorem}{Theorem}[section]
\newtheorem{proposition}[theorem]{Proposition}
\newtheorem{lemma}[theorem]{Lemma}

\newtheorem{question}[theorem]{Question}

\theoremstyle{definition}

\newtheorem{definition}[theorem]{Definition}

\newtheorem{remark}[theorem]{Remark}

\def\({\left(}
\def\){\right)}

\def\lp{\left(}
\def\rp{\right)}

\def\({\left(}
\def\){\right)}

\def\lp{\left(}
\def\rp{\right)}

\newcommand{\Log}{\mathrm{Log}}

\title[Distribution of Alternating Sums]{Distribution of Alternating Sums of Parts in Partitions}
\date{}

\author{William Craig}
\address{US Naval Academy}
\email{wcraig@usna.edu}

\author{Runqiao Li}
\address{Penn State University}
\email{rml5856@psu.edu or runqiaoli@outlook.com}

\begin{document}

\begin{abstract}
    Recently, many authors have investigated how various partition statistics distribute as the size of the partition grows. In this work, we look at a particular statistic arising from the recent rejuvenation of MacMahon's partition analysis. More specifically, we compute all the moments of the alternating sum statistic for partitions. We prove this results using the Circle Method. We also propose a general framework for studying further questions of this type that may avoid some of the complications that arise in traditional approaches to the distributions of partition statistics, and we comment on the utility, comparative ease and opportunities to generalize to very broad settings.
\end{abstract}

\dedicatory{Commemorating the 85th birthdays of George E. Andrews and Bruce C. Berndt}

\maketitle

\section{Introduction}

A {\it partition} of a non-negative integer $n$ is a non-increasing list $\lp \lambda_1, \lambda_2, \dots, \lambda_\ell \rp$ of positive integers such that $|\lambda| := \lambda_1 + \lambda_2 + \dots + \lambda_\ell = n$. We will let $\ell = \ell(\lambda)$ be the number of parts of $\lambda$, $\lambda \vdash n$ denote that $\lambda$ is a partition of $n$, and $\mathcal P$ be the collection of all partitions. The statistical theory of partitions has been a very fruitful area of analytic combinatorics, and indeed is often the first case of much larger theories in all aspects of combinatorics. Many interesting results have been proven on the limiting shapes of Young diagrams for partitions, including works by Szalay and Tur\'{a}n \cite{SzalayTuran1,SzalayTuran2,SzalayTuran3}, Temperely \cite{Temperely}, Vershik \cite{Vershik}, and Bogachev \cite{Bogachev}. In particular, works beginning with Erd\H{o}s and Lehner \cite{ErdosLehner} and greatly strengthened by the framework of Fristedt \cite{Fristedt} have given much attention to how $\ell(\lambda)$ distributes among partitions of $n$. Such ideas have been extended to a number of other combinatorial settings; see for instance the work of Bridges and Bringmann \cite{BridgesBringmann} applying Fristedt's ideas to unimodal sequences.

Our objective is to study these questions for new statistics on partitions, motivated by studies of distributions on Betti numbers for Hilbert schemes \cite{GriffinOnoRolenTsai} and of hook lengths in partitions \cite{GriffinOnoTsai} and restricted classes of partitions \cite{BallantineBursonCraigFolsomWen,CraigDawseyHan,CraigOnoSingh}. We will in particular consider this question for the {\it alternating sum}
\begin{align*}
    a(\lambda) := \lambda_1 - \lambda_2 + \lambda_3 - \dots = \sum_{j \geq 1} (-1)^j \lambda_j.
\end{align*}
Note that we let $\lambda_j = 0$ for $j > \ell(\lambda)$. This statistic, which is equivalent to the number of odd parts in the conjugate of $\lambda$, was studied in \cite{Andrews2,BerkovichUncu,KangLiWang,WangZheng1} and has emerged out of a stream of work on Schmidt-type statistics for partitions \cite{AndrewsKeith, Li,WangZheng2} that has arisen from recent developments in MacMahon's partition analysis \cite{AndrewsPaule2,AndrewsPaule3,AndrewsPaule4}. Our main result on $a(\lambda)$ is phrased in terms of the {\it moments} of $a(\lambda)$:

\begin{definition}
    We define for any integers $m \geq 1$, $n \geq 0$ the sequences
    \begin{align*}
        A_m(n) := \sum_{\lambda \vdash n} a(\lambda)^m
    \end{align*}
    We also define the moments of $a(\lambda)$, which are given by
    \begin{align*}
        \mathbb{E}_m(n) := \dfrac{A_m(n)}{p(n)}.
    \end{align*}
\end{definition}

Our main results then describes an asymptotic formula for $\mathbb{E}_m(n)$ below.

\begin{theorem} \label{T: Main Theorem}
    Let $m \geq 1$. Then we have, as $n \to \infty$, that 
    \begin{align*}
        \mathbb{E}_m(n) \sim \dfrac{6^{\frac{m}{2}}}{2^m \pi^m} n^{\frac{m}{2}} \log^m\lp \dfrac{\sqrt{6n}}{\pi} \rp.
    \end{align*}
\end{theorem}

\begin{remark}
   By using an alternative interpretation of the generating function for $a(\lambda)$, we also derive the moments for the number of odd parts in partitions. The distribution for this statistic is proven using different techniques by Yakubovich \cite[Theorem 8.1]{Yakubovich} and probably could again be proven using sieving techniques like those in \cite{ErdosLehner,GriffinOnoRolenTsai}. Our technique is useful, as we discuss later, in its wide ability to be generalized and comparatively little theory required. In particular, the aforementioned result of Yakubovich states that
   \begin{align*}
        \lim\limits_{n \to \infty} \#\left\{ \lambda \vdash n : \dfrac{\pi}{\sqrt{6n}} a(\lambda) - \dfrac{1}{4} \log(n) + \dfrac{1}{2} \log\lp \dfrac{2\pi}{\sqrt{6}} \rp \leq x \right\} = \mathrm{Erfc}\lp e^{-x} \rp,
    \end{align*}
    where
    \begin{align*}
        \mathrm{Erfc}(u) = \dfrac{2}{\sqrt{\pi}} \int_u^\infty e^{-z^2} dz.
    \end{align*}
   By determining additional asymptotic terms, which could be done in principle with our method, this same distribution can be derived. We note also the work of Goh and Schmutz on the distribution of the number of distinct part sizes in partitions \cite{GohSchmutz}, which gives a normal distribution. An visual example of how this distribution plays out is given below.
\end{remark}

\begin{figure}[h] 
\begin{center}
\includegraphics[width=10cm,height=5cm]{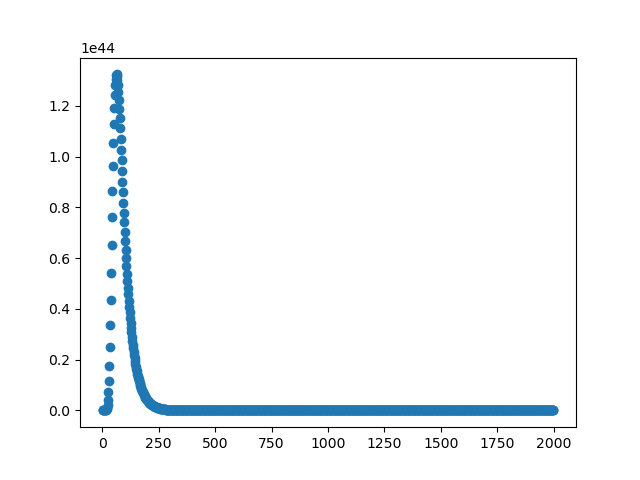} 
\end{center}
\caption{Distribution of all the values of $a(\lambda)$ among the partitions of 2000.}
\label{Figure 1}
\end{figure}

The remainder of the paper is organized as follows. In Section \ref{S: Preliminaries}, we provide preliminary details needed for this work, particularly the construction of generating function for $a(\lambda)$ and asymptotic results derived from Euler--Maclaurin summation. In Section \ref{S: Asymptotic} we apply results of Section \ref{S: Preliminaries} to the specific context needed for the paper. In Section \ref{S: Proof}, we prove Theorem \ref{T: Main Theorem}. Finally, in Section \ref{S: Further Discussion} we lay out an argument in favor of our approach to distributions for partition statistics, comparing our result with results of Bringmann, Mahlburg and Rhoades \cite{BringmannMahlburgRhoades} on ranks and cranks of partitions, and arguing for the opportunities that come out of uniting the combinatorial study of MacMahon's partition analysis with the analytic technology of the Circle Method.

\section*{Acknowledgments}

The authors thank Ken Ono for introducing the authors and for discussions which led to this project and we thank Daniel Parry for helpful comments and additional references which corrected typos and improved the exposition of the paper, and we thank Walter Bridges for spotting a significant error in the original manuscript. The second author would also like to thank Amita Malik for her very productive discussion. The views expressed in this article are those of the authors and do not reflect the official policy or position of the U.S. Naval Academy, Department of the Navy, the Department of Defense, or the U.S. Government.

\section{Preliminaries} \label{S: Preliminaries}

In this section, we provide various lemmas which are needed for the main proof; in particular, we focus on the combinatorial construction of a generating function and the asymptotic analysis of this generating function via Euler--Maclaurin summation.

\subsection{Generating function}

The construction of the generating function for the alternating part sum statistic relies centrally on the theory of partition analysis introduced by MacMahon in his seminal volumes on combinatorial analysis \cite{MacMahon}, which has been developed to a very high degree by Andrews, Paule, and their collaborators in recent decades. In particular, the theory hinges on the {\it Omega operator} defined by
\begin{align*}
    \Omega_{\geq} \sum_{s_1 = -\infty}^\infty \cdots \sum_{s_r = -\infty}^\infty A_{s_1,\dots,s_r} \lambda_1^{s_1} \cdots \lambda_r^{s_r} = \sum_{s_1 = 0}^\infty \cdots \sum_{s_r = 0}^\infty A_{s_1,\dots,s_r},
\end{align*}
where the domain of the $A_{s_1,\dots,s_r}$ is the field of rational functions over $\C$ in several complex variables and the $\lambda_i$ are restricted to annuli of the form $1 - \varepsilon < |\lambda_i| < 1 + \varepsilon$. The basic purpose of the MacMahon operator, from the perspective of combinatorial analysis and partition theory, is to create generating functions for interesting statistics by encoding restrictions on the components of combinatorial objects within the exponents of the parameters $\lambda_i$. The construction of a generating function tracking the statistic $a(\lambda)$ (and indeed any number of related statistics which might be studied from a combinatorial lens) relies crucially on this framework, which we now lay out.

\begin{lemma}\label{elimination}
For any integer $A\geq0$, we have
\begin{align*}
\underset{\geq}{\Omega}\frac{\lambda^{-A}}{(1-x\lambda)(1-\frac{y}{\lambda})}=\frac{x^{A}}{(1-x)(1-xy)}.
\end{align*}
\end{lemma}

\begin{proof}
By the definition of the Omega operator, we have
\begin{align*}
\underset{\geq}{\Omega}\frac{\lambda^{-A}}{(1-x\lambda)(1-\frac{y}{\lambda})}=&\underset{\geq}{\Omega}\sum_{n=0}^{\infty}\sum_{m=0}^{\infty}x^{n}y^{m}\lambda^{n-m-A}=\sum_{m=0}^{\infty}\sum_{n=m+A}^{\infty}x^{n}y^{m}=\frac{x^{A}}{(1-x)(1-xy)}.
\end{align*}
This completes the proof.
\end{proof}

Our techniques begin with a two-variable generating function that tracks simultaneously $a(\lambda)$ and $|\lambda|$, which is critical to our method.

\begin{lemma}\label{generating function}
    We have
    \begin{align*}
        P(z;q) := \sum_{\lambda \in \mathcal P} z^{a(\lambda)} q^{|\lambda|} = \dfrac{1}{\lp zq;q^2 \rp_\infty \lp q^2;q^2 \rp_\infty}.
    \end{align*}
\end{lemma}

\begin{remark}
    Note this is also the generating function giving the distribution of the number of odd part sizes appearing in partitions. Thus, the distribution determined by our work apply also to this statistic.
\end{remark}

\begin{proof}
    We start by consider the generating function for partitions with length bounded by $n$. For a partition $\lambda=(\lambda_1,\lambda_2,\ldots,\lambda_n)$, define the weight $x^{\lambda}:=x_1^{\lambda_1}x_2^{\lambda_2}\cdots x^{\lambda_n}$. We have
    \begin{align*}
    \sum_{\substack{\lambda\in\mathcal{P}\\\ell(\lambda)\leq n}}x^{\lambda}=&\sum_{a_1\geq a_2\geq\cdots\geq a_n\geq0}x_1^{a_1}x_2^{a_2}\cdots x_n^{a_n}\\
    =&\underset{\geq}{\Omega}\sum_{a_1,a_2,\ldots,a_n\geq0}x_1^{a_1}x_2^{a_2}\cdots x_n^{a_n}\lambda_1^{a_1-a_2}\lambda_2^{a_2-a_3}\cdots\lambda_{n-1}^{a_{n-1}-a_n}\lambda_n^{a_n}\\
    =&\underset{\geq}{\Omega}\frac{1}{(1-x_1\lambda_1)(1-\frac{x_2\lambda_2}{\lambda_1})\cdots(1-\frac{x_{n-1}\lambda_{n-1}}{\lambda_{n-2}})(1-\frac{x_n\lambda_n}{\lambda_{n-1}})},
    \end{align*}
    which is usually called the ``crude form" of the generating function. Our next step is to eliminate the $\lambda_j$'s one by one, starting with $\lambda_1$. Applying Lemma \ref{elimination} with $A=0$ to eliminate first $\lambda_1$, then $\lambda_2$, and so on proceeding to $\lambda_n$, we have
    \begin{align*}
     \sum_{\substack{\lambda\in\mathcal{P}\\\ell(\lambda)\leq n}}x^{\lambda}=&\underset{\geq}{\Omega}\frac{1}{(1-x_1\lambda_1)(1-\frac{x_2\lambda_2}{\lambda_1})\cdots(1-\frac{x_{n-1}\lambda_{n-1}}{\lambda_{n-2}})(1-\frac{x_n\lambda_n}{\lambda_{n-1}})}\\
     =&\frac{1}{(1-x_1)} \cdot \underset{\geq}{\Omega}\frac{1}{(1-{x_1x_2\lambda_2})(1-\frac{x_3\lambda_3}{\lambda_2})\cdots(1-\frac{x_{n-1}\lambda_{n-1}}{\lambda_{n-2}})(1-\frac{x_n\lambda_n}{\lambda_{n-1}})}\\
     =&\frac{1}{(1-x_1)(1-x_1x_2)} \cdot \underset{\geq}{\Omega}\frac{1}{(1-{x_1x_2x_3\lambda_3})\cdots(1-\frac{x_{n-1}\lambda_{n-1}}{\lambda_{n-2}})(1-\frac{x_n\lambda_n}{\lambda_{n-1}})}\\
     &\vdots\\
     =&\frac{1}{(1-x_1)(1-x_1x_2)\cdots(1-x_1\cdots x_n)},
    \end{align*}
    Let $n\to\infty$ and define $X_i=x_1x_2\cdots x_i$ for $i\geq1$, we have
    $$\sum_{\lambda\in\mathcal{P}}x^{\lambda}=\prod_{n=1}^{\infty}\frac{1}{1-X_n}.$$
    Let $x_{2i-1}\to zq$ and $x_{2i}\to q/z$ for all $i\geq1$, which implies $X_{2i-1}=zq^{2i-1}$ and $X_{2i}=q^{2i}$, we have
    $$\sum_{\lambda\in\mathcal{P}}z^{a(\lambda)}q^{|\lambda|}=\frac{1}{(zq;q^2)_{\infty}(q^{2};q^{2})_{\infty}}.$$
    This completes the proof.
\end{proof}

\begin{remark}
There are other ways to get this two-variable generating function. For instance, one can easily derive it from the main result of \cite{Boulet}. However, here we present a complete proof by MacMahon's partition analysis to serve as an example, so that the reader can study such generating function for different type of partitions by following the similar steps.
\end{remark}

\subsection{Dilogarithm function}

We recall the {\it dilogarithm function} $\mathrm{Li}_2\lp z \rp$, which is defined by
\begin{align*}
    \mathrm{Li}_2(z) := \sum_{k \geq 1} \dfrac{z^k}{k^2}
\end{align*}
for $|z|<1$ and is extended by analytic continuation to $\mathbb{C} \backslash [1,\infty)$. This function is among the most important non-elementary functions, and plays a major role in many partition asymptotic problems. We require for our purposes the distribution property \cite[pg. 9]{ZagierDilog}
\begin{align} \label{eq: distribution property}
	\mathrm{Li}_2\lp x \rp = n \sum_{z^n = x} \mathrm{Li}_2\lp z \rp.
\end{align}

\subsection{The saddle point method}

We recall here a classical statement of the saddle point method (or Laplace's method) which is used for computing asymptotics of exponential integrals with oscillation.

\begin{theorem}[{\cite[Section 1.1.5]{Pinsky}}] \label{T: Saddle point method}
    Let $A,B : [a,b] \to \C$ be continuous functions. Suppose there is $x_0 \in [a,b]$ such that for $x \not = x_0$ in $[a,b]$, we have $\mathrm{Re}(x) < \mathrm{Re}(x_0)$, and that
    \begin{align*}
        \lim_{x \to x_0} \dfrac{B(x) - B(x_0)}{(x-x_0)^2} = -k \in \C
    \end{align*}
    where $\mathrm{Re}(k)>0$. Then as $t \to \infty$, we have
    \begin{align*}
        \int_a^b A(x) e^{tB(x)} dx = e^{tB(x_0)}\lp A(x_0) \sqrt{\dfrac{\pi}{tk}} + o\lp \dfrac{1}{\sqrt{t}} \rp \rp.
    \end{align*}
\end{theorem}

\subsection{Euler--Maclaurin summation}

To state these results, we will say that a function $f(w)$ satisfies the asymptotic $f(w) \sim \sum_{n \geq n_0} c_n w^n$ as $w \to 0$ if $f(w) - \sum_{n = n_0}^{N - 1 + n_0} c_n w^n = O_N\lp w^N \rp$ for each $N \geq 1$. We also define the integrals (for $A \in \R)$
\begin{align*}
    I_f := \int_0^\infty f(x) dx \ \ \ \text{and} \ \ \ I_{f,A}^*:=\int_{0}^{\infty} \left(f(u)-\sum_{n=n_0}^{-2}c_{n}u^n-\frac{c_{-1}e^{-Au}}{u}\right)du.
\end{align*}
where convergent. We define the regions
\begin{align*}
    D_\theta := \{ re^{i\alpha} : r > 0 \mbox{ and } |\alpha|\le \theta  \}.
\end{align*}
We begin with an exact version of Euler--Maclaurin summation as stated by Bringmann, Jennings--Shaffer and Mahlburg \cite{BringmannJenningsShafferMahlburg}.

\begin{lemma}[{\cite[Equation 5.8]{BringmannJenningsShafferMahlburg}}] \label{L: EM Exact}
	Suppose that $0 \leq \theta < \frac{\pi}{2}$. Let $f : \C \to \C$ be holomorphic in a domain containing $D_\theta$, so that in particular $f$ is holomorphic at the origin, and assume that f and all of its derivatives are of sufficient decay (i.e. decays at least as quickly as $|w|^{-1-\varepsilon}$ for some $\varepsilon>0$). Then for $a \in \R$ and $N \in \N_0$, we have
	\begin{multline*}
		\sum_{m\geq0} f((m+a)w) = \frac{I_f}{w} - \sum_{n=0}^{N-1} \frac{B_{n+1}(a)f^{(n)}(0)}{(n+1)!} w^n - \sum_{k\ge N} \frac{f^{(k)}(0)a^{k+1}}{(k+1)!}w^k\\
		- \frac{w^N}{2\pi i} \sum_{n=0}^{N-1} \frac{B_{n+1}(0)a^{N-n}}{(n+1)!} \int_{C_R(0)} \frac{f^{(n)}(z)}{z^{N-n}(z-aw)} dz - (-w)^{N-1} \hspace{-0.1cm} \int_{aw}^{w\infty} \frac{f^{(N)}(z)\widetilde B_N\left(\frac zw-a\right)}{N!} dz,
	\end{multline*}
	where $\widetilde B_n(x):=B_n(x-\lfloor x\rfloor)$ and $C_R(0)$ denotes the circle of radius $R$ centered at the origin, where $R$ is such that $f$ is holomorphic in $C_{R}(0)$.
\end{lemma}

\begin{remark}
    Note that none of the results in this section (apart from where explicitly stated) require local convergence of a Taylor/Laurent series for $f$, only a valid asymptotic expansion near zero.
\end{remark}

This exact formula serves as the basis for computing useful asymptotic expansions. Although we will once use Lemma \ref{L: EM Exact} in the paper, it is also greatly useful in deriving simpler asymptotic formulas. We state two such lemmas below. The first is essentially a restatement of \cite[Lemma 2.2]{BringmannCraigMalesOno} with additional comments on derivatives, and the second a special case of this.

\begin{lemma} \label{L: EM with Principal Parts}
	Assume the notation of Lemma \ref{L: EM Exact}, and that $f$ is meromorphic at the origin and that $f$ and all its derivatives are of sufficient decay. Then we have
	\begin{align*}
        \sum_{n=0}^\infty f((n+a)w) &\sim \sum_{n=n_0}^{-2} c_{n} \zeta(-n,a)w^{n} + \frac{I_{f,A}^*}{w} \\ &-\frac{c_{-1}}{w} \left( \Log \left(Aw \right) +\psi(a)+\gamma \right)-\sum_{n=0}^\infty c_n \frac{B_{n+1}(a)}{n+1} w^n,
	\end{align*}
	as $w \rightarrow 0$ uniformly in $D_\theta$. Furthermore, if $j \geq 1$ is given and if the $j$th derivative of the series above has an asymptotic expansion, then
    \begin{align*}
        \dfrac{d^j}{dw^j} \sum_{n=0}^\infty f((n+a)w) &\sim \dfrac{d^j}{dw^j} \bigg[ \sum_{n=n_0}^{-2} c_{n} \zeta(-n,a)w^{n} + \frac{I_{f,A}^*}{w} \\ &-\frac{c_{-1}}{w} \left( \Log \left(Aw \right) +\psi(a)+\gamma \right)-\sum_{n=0}^\infty c_n \frac{B_{n+1}(a)}{n+1} w^n \bigg],
	\end{align*}
\end{lemma}

\begin{proof}
    This is \cite[Lemma 2.2]{BringmannCraigMalesOno} along with the fact that term-by-term integration of asymptotic expansions is valid.
\end{proof}

The fact that the exchange of derivatives is permitted is not difficult to establish in many cases, which leads to the following result.

\begin{lemma}[{\cite[Corollary 2.3]{CraigOnoSingh}}] \label{L: EM with derivatives}
    Adopt the notation of Lemma \ref{L: EM with Principal Parts}, and suppose $f(w)$ is holomorphic in a domain containing $D_\theta$, and that $f$ and all its derivatives are of sufficient decay as $|w| \to \infty$ (i.e. decays at least as quickly as $|w|^{-1-\varepsilon}$ for some $\varepsilon>0$). Suppose also that $f$ has a Taylor expansion $f(w) = \sum_{n \geq 0} c_n w^n$ near $w=0$. Then for $a \in \mathbb{R}_{>0}$, we have as $w \to 0$ in $D_\theta$ that
    \begin{align*}
        \sum_{m \geq 0} f\lp (m+a)w \rp \sim \dfrac{I_f}{w} - \sum_{n \geq 0} c_n \dfrac{B_{n+1}(a)}{n+1} w^n,
    \end{align*}
    where $B_n(x)$ are the Bernoulli polynomials. Furthermore, for any $j \geq 0$ we have as $w \to 0$ in $D_\theta$ that
    \begin{align*}
		\dfrac{d^j}{dw^j} \sum_{m\geq0} f((m+a)w) \sim \frac{(-1)^j j!}{w^{j+1}} I_f - \sum_{n \geq 0} c_n \dfrac{n!}{(n-j)!} \cdot \dfrac{B_{n+1}(a)}{n+1} w^n.
	\end{align*}
\end{lemma}

\begin{proof}
    The first statement is a corollary immediately of Lemma \ref{L: EM with Principal Parts}. Observing that the proper decay conditions are in place, we exchange derivatives and sums and we can write
     \begin{align*}
        \dfrac{d^j}{dw^j} \sum_{m \geq 0} f\lp (m+a)w \rp = \sum_{m \geq 0} \lp m+a \rp^j f^{(j)}\lp (m+a)w \rp &= \dfrac{1}{w^j} \sum_{m \geq 0} g\lp (m+a)w \rp,
    \end{align*}
    where we let $g(w) := w^j f^{(j)}(w)$. Note then from the Taylor expansion of $f$ as $w \to 0$ that 
    \begin{align*}
        g(w) = w^j \dfrac{d^j}{dw^j} \sum_{n \geq 0} c_n w^n = \sum_{n \geq j} \dfrac{n!}{(n-j)!} c_n w^n =: \sum_{n \geq 0} d_n w^n,
    \end{align*}
    where $d_0 = \cdots = d_{j-1} = 0$ and $d_n = \frac{n!}{(n-j)!} c_n$ for $n \geq j$. Then by application of the $j=0$ case of the lemma, we obtain
    \begin{align*}
        \sum_{m \geq 0} g\lp (m+a) w \rp \sim \dfrac{I_g}{w} + \sum_{n \geq j} c_n \dfrac{n!}{(n-j)!} \cdot \dfrac{B_{n+1}(a)}{n+1} w^n
    \end{align*}
    as $w \to 0$ in $D_\theta$. Evaluation of $I_g$ in terms of $I_f$ is simple using integration by parts, and the proof follows.
\end{proof}

\subsection{Applications of Euler--Maclaurin summation}

We use these results to compute asymptotics for infinite products built from pieces of the shape $\lp z q^a; q^b \rp_\infty^{-1}$. Results of this type are the main technical input to the saddle point method approximation for the generating functions of partition statistics. More specifically, we compute asymptotics for these products as $q$ approaches a root of unity, which we denote using the standard notation $\zeta_k^h := e^{\frac{2\pi i h}{k}}$. We will first handle cases which do not involve $z$, which are somewhat simpler.

\begin{lemma} \label{L: No z}
    Let $q = e^{-w}$ be in the unit disk, $0 \leq \frac hk < 1$ be a fraction in simplified terms, and let $1 \leq a \leq b$ be integers with $b$ coprime to $k$. Then as $q \to \zeta_k^h$ in the region $\zeta_k^h \cdot D_\theta$, we have
    \begin{align*}
        \Log\lp q^a; q^b \rp_\infty^{-1} \begin{cases}
            \sim \dfrac{\zeta(2)}{bw} + b I_{f_{a,b},A}^* - \lp \dfrac{b}{2} - a \rp \lp \Log\lp Aw \rp + \psi\lp \dfrac ab \rp + \gamma \rp + O(w) & \text{if} \ \ \dfrac hk = 0, \\[+0.2cm]
            = o\lp w^{-1} \rp & \text{if} \ \ \dfrac hk \not = 0,
        \end{cases}
    \end{align*}
    where $f_{a,b}(w) := \frac{e^{-aw/b}}{w\lp 1 - e^{-w} \rp}$ and $A>0$ is freely chosen.
\end{lemma}

\begin{remark}
    We note that the integrals $I_{f_{a,b},A}$ might be evaluated using techniques like those in \cite[Lemma 2.3]{BringmannCraigMalesOno}. We do not require this level of generality, but do apply this result for one particular evaluation later in this work. The true order of $o(w^{-1})$ would be of size $\Log(w)$ or of constant size depending on asymptotic expansions of underlying functions.
\end{remark}

\begin{proof}
    We start with $\frac hk = \frac 01$, where by expanding $\Log\lp 1 - q^{2j} \rp$ as a series we obtain
    \begin{align*}
        \Log\lp q^a; q^b \rp_\infty^{-1} = \sum_{j \geq 0} \sum_{m \geq 1} \dfrac{(q^{a+bj})^m}{m} = \sum_{m \geq 1} \dfrac{q^{am}}{m\lp 1 - q^{bm} \rp} = bw \sum_{m \geq 1} \dfrac{e^{-amw}}{bmw\lp 1 - e^{-bmw} \rp}.
    \end{align*}
    Now, applying Lemma \ref{L: EM with Principal Parts} in the case $f_{a,b}(w) := \frac{e^{- aw/b}}{w\lp 1 - e^{-w} \rp} = \frac{1}{w^2} + O\lp w^{-1} \rp$, we obtain
    \begin{align*}
        \Log\lp q^a; q^b \rp_\infty^{-1} = bw\lp \dfrac{\zeta(2)}{(bw)^2} + \dfrac{I_{f_{a,b},A}^*}{w} - \dfrac{\frac{1}{2} - \frac{a}{b}}{w}\lp \Log\lp Aw \rp + \psi(a) - \gamma \rp + O(1) \rp
    \end{align*}
    as $w \to 0$ in $D_\theta$. This completes the proof for $q \to 1$.
    
    If on the other hand we pick $\zeta_k^h \not = 1$, we change variables $q \mapsto \zeta_k^h q$, so that now we may take $q \to 1$ in $D_\theta$ instead of $q \to \zeta_k^h$ in $\zeta_k^h \cdot D_\theta$. Doing this, we observe that
    \begin{align*}
        \Log\lp \lp \zeta_k^h q \rp^a; \lp \zeta_k^h q \rp^b \rp_\infty^{-1} &= \sum_{j \geq 0} \sum_{m \geq 1} \dfrac{\lp \zeta_k^{ah + bhj} q^{a+bj} \rp^m}{m} = \sum_{m \geq 1} \dfrac{\zeta_k^{ahm} q^{am}}{m\lp 1 - \zeta_k^{bh} q^{bm} \rp} \\ &= bw \sum_{r=1}^k \sum_{\ell \geq 0} f_{a,b,r,h,k}\lp \lp \ell + \frac rk \rp bkw \rp,
    \end{align*}
    where
    \begin{align*}
        f_{a,b,r,h,k}(w) := \dfrac{\zeta_k^{ahr} e^{-\frac ab w}}{w \lp 1 - \zeta_k^{bh} e^{-w} \rp} = \dfrac{\zeta_k^{ahr}}{\lp \zeta_k^{bh} - 1 \rp w} + O(1)
    \end{align*}
    since $bh$ is necessarily coprime to $k$. It then follows from Lemma \ref{L: EM with Principal Parts} that
    \begin{align*}
        \Log\lp \lp \zeta_k^h q \rp^a; \lp \zeta_k^h q \rp^b \rp_\infty^{-1} = o\lp w^{-1} \rp
    \end{align*}
    as required.
\end{proof}

For the sake of completeness, i.e. to include the most general asymptotic results which are relevant to the genre of literature to which we contribute, we include the following additional lemma which may be very useful to those using a saddle point approach to distributions, as we see in \cite{CraigOnoSingh,GriffinOnoTsai,GriffinOnoRolenTsai}, and others.

\begin{lemma} \label{L: One z}
    Let $z \in \mathbb{C}$ be nonzero and $q = e^{-w}$ be in the unit disk. Let $0 \leq \frac hk < 1$ be a fraction in simplified terms, and let $1 \leq a \leq b$ be integers with $b$ coprime to $k$. Then as $q \to \zeta_k^h$ in the region $\zeta_k^h \cdot D_\theta$ for some $\theta > 0$, we have as long as $z^k \not \in [1,\infty)$ that
    \begin{align*}
        \Log\lp z q^a; q^b \rp_\infty^{-1} = \dfrac{\mathrm{Li}_2\lp z^k \rp}{bk^2w} + O_\theta(k).
    \end{align*}
\end{lemma}

\begin{remark}
    One can also compute such asymptotics in cases where $\gcd(b,k) > 1$ with a little more work, this can be found for instance in \cite{Parry}.
\end{remark}

\begin{proof}
    By the proof of Lemma \ref{L: EM with derivatives}, we need only consider the case without derivatives. Note that taking $q \to \zeta_k^h$ in $\zeta_k^h \cdot D_\theta$ is equivalent to taking $q \to 1$ in $D_\theta$ after replacing $q \mapsto \zeta_k^h q$ in generating functions. We therefore consider asymptotics for function $\lp z \zeta_k^{ah} q^a; \zeta_k^{bh} q^b \rp_\infty^{-1}$ as $q \to 1$ in $D_\theta$ using Lemma \ref{L: EM with derivatives}. To do so, we first note
    \begin{align*}
        \Log\lp z \zeta_k^{ah} q^a; \zeta_k^{bh} q^b \rp_\infty^{-1} = - \sum_{m \geq 0} \Log\lp 1 - z \zeta_k^{h(a+bm)} q^{a+bm} \rp.
    \end{align*}
    Now, we decompose $m = k\ell + c$, $\ell \geq 0$ and $0 \leq c \leq k-1$, so we may then write
    \begin{align*}
        \Log\lp z \zeta_k^{ah} q^a; \zeta_k^{bh} q^b \rp_\infty^{-1} &= - \sum_{\ell \geq 0} \sum_{c=0}^{k-1} \Log\lp 1 - z \zeta_k^{ah + bh(k\ell + c)} q^{a+b(k\ell+c)} \rp \\ &= - \sum_{c=0}^{k-1} \sum_{\ell \geq 0} \Log\lp 1 - z \zeta_k^{h(a + bc)} e^{-\lp \frac{a+bc}{bk}+\ell \rp bkw} \rp.
    \end{align*}
    Now, let
    \begin{align*}
        f_{a,b,c,h,k,z}\lp w \rp := - \Log\lp 1 - z \zeta_k^{h(a+bc)} e^{-w} \rp,
    \end{align*}
    so that we have
    \begin{align*}
        \Log\lp z \zeta_k^{ah} q^a; \zeta_k^{bh} q^b \rp_\infty^{-1} = \sum_{c=0}^{k-1} \sum_{\ell \geq 0} f_{a,b,c,h,k,z}\lp \lp \ell + \frac{a+bc}{bk} \rp bkw \rp.
    \end{align*}
    Since $z^k \not = 1$, we have $f_{a,b,c,h,k,z}\lp w \rp$ is holomorphic at $w=0$, and therefore by Lemma \ref{L: EM with derivatives} we obtain
    for each $j \geq 0$ that
    \begin{align*}
        \sum_{\ell \geq 0} f_{a,b,c,h,k,z}\lp \lp \ell + \frac{a+bc}{bk} \rp bkw \rp = \dfrac{1}{bkw} \int_0^\infty f_{a,b,c,h,k,z}(x) dx + O_\theta(1).
    \end{align*}
    Using the integral evaluation
    \begin{align*}
        - \int_0^\infty \Log\lp 1 - z \zeta_k^h e^{-x} \rp dx = \mathrm{Li}_2\lp z \zeta_k^h \rp,
    \end{align*}
    which is valid since $z \zeta_k^h \not \in [1,\infty)$, we have
    \begin{align*}
        \Log\lp z \zeta_k^{ah} q^a; \zeta_k^{bh} q^b \rp_\infty^{-1} = \sum_{c = 0}^{k-1} \dfrac{\mathrm{Li}_2\lp z \zeta_k^{h(a+bc)} \rp}{bkw} + O_\theta(k).
    \end{align*}
    Since we assume that $b$ is coprime to $k$, $h(a+bc)$ runs over all congruence classes modulo $k$ exactly once as we run $0 \leq c \leq k-1$, and therefore by the distribution property \eqref{eq: distribution property} we obtain
    \begin{align*}
        \Log\lp z \zeta_k^{ah} q^a; \zeta_k^{bh} q^b \rp_\infty^{-1} = \dfrac{\mathrm{Li}_2\lp z^k \rp}{bk^2 w} + O_\theta(k),
    \end{align*}
    which completes the proof.
\end{proof}

\section{Asymptotics for Generating Functions} \label{S: Asymptotic}

In this section, we apply asymptotic results from Section \ref{S: Preliminaries} to the generating functions with which we are concerned in order to prove Theorem \ref{T: Main Theorem}.

\begin{lemma}
    For any $m \geq 1$, we have the generating function identity
    \begin{align*}
        \mathcal A_m(q) := \sum_{n \geq 0} A_m(n) q^n = \lp z\dfrac{\partial}{\partial z} \rp^m \bigg|_{z=1} P(z;q).
    \end{align*}
\end{lemma}

\begin{proof}
    This is a straightforward calculation with derivatives.
\end{proof}

\begin{proposition} \label{P: Major arc condition}
    Let $m \geq 1$. As $|w| \to 0$ in $D_\delta$ for $0 < \delta < \frac{\pi}{2}$, we have
    \begin{align*}
        \mathcal A_m\lp q \rp = \dfrac{\Log^m(1/w)}{2^m \sqrt{2\pi} w^{m - \frac 12}} e^{\frac{\pi^2}{6w}} \lp 1 + O\lp \dfrac{1}{\Log(1/w)} \rp \rp.
    \end{align*}
\end{proposition}

\begin{proof}
    It is straightforward to see in the first place that
    \begin{align*}
        z \dfrac{\partial}{\partial z} P(z;q) = G(z;q) P(z;q),
    \end{align*}
    where
    \begin{align*}
        G(z;q) := \sum_{j \geq 0} \dfrac{zq^{2j+1}}{1 - zq^{2j+1}}.
    \end{align*}
    By induction of $m$, it is clear that there is a sequence of polynomials $F_m(x_0,\dots,x_{m-1})$ such that
    \begin{align*}
        \lp z \dfrac{\partial}{\partial z} \rp^m P(z;q) = F_m\lp G(z;q), G^{(1)}(z;q), \dots, G^{(m-1)}(z;q) \rp \cdot P(z;q),
    \end{align*}
    where here and throughout derivatives signify derivatives in the $z$ variable. In particular, operating under the induction hypothesis for a particular $m \geq 1$, we have
    \begin{align*}
        &\lp z \dfrac{\partial}{\partial z} \rp^{m+1} P(z;q) = z \dfrac{\partial}{\partial z} \left[ F_m\lp G(z;q), \dots, G^{(m-1)}(z;q) \rp \cdot P(z;q) \right] \\ &= \lp z \dfrac{\partial}{\partial z} \left[ F_m\lp G(z;q), \dots, G^{(m-1)}(z;q) \rp \right] + F_m\lp G(z;q), \dots, G^{(m-1)}(z;q) \rp G(z;q)\rp P(z;q).
    \end{align*}

    It is also straightforward to demonstrate by induction on $k$ that for each $k \geq 1$,
    \begin{align*}
        \left.\lp z\frac{\partial}{\partial z}\rp^kG(z;q)\right|_{z=1}=\sum_{j\geq0}\frac{E_{k}(q^{2j+1})}{(1-q^{2j+1})^{k+1}},
    \end{align*}
    where $E_k(q)$ are the Eulerian polynomials defined by
    \begin{align*}
        \sum_{j \geq 1} j^k x^j = \dfrac{E_k(x)}{\lp 1 - x \rp^{k+1}}.
    \end{align*}
    
    Now,for $q = e^{-w}$, it is straightforward to prove that
    \begin{align*}
        \dfrac{E_k\lp q^{2j+1} \rp}{\lp 1 - q^{2j+1} \rp^{k+1}} = \dfrac{k!}{2^{k+1} w^{k+1}}\lp 1 + O(w) \rp,
    \end{align*}
    and therefore by applying Lemma \ref{L: EM with Principal Parts}, we then see that for $q = e^{-w}$, then as $q \to 1$ on the major arc path we have (using known identities for the digamma function)
    \begin{align*}
        G^{(k)}(1;q) = \begin{cases}
            \dfrac{k! \ \zeta(k+1)}{w^{k+1}} \lp 1 + O(w) \rp & \textrm{for } k \geq 1, \\[+0.5cm]
            \dfrac{I_{f,A}^*}{w} - \dfrac{1}{2w}\lp \Log\lp Aw \rp - \log(4) \rp, & \textrm{for } k = 0,
        \end{cases}
    \end{align*}
    for any $A > 0$ and where $f(w) := \frac{e^{-2w}}{1 - e^{-2w}}$. Choosing $A=1$, we have
    \begin{align*}
        G(1;q) = \dfrac{\Log(1/w)}{2w} + O\lp \dfrac{1}{w} \rp
    \end{align*}

    Now, as a sample, the first few polynomials $F_m$ are given by
    \begin{align*}
        F_1\lp G(1;q) \rp &= G(1;q) = \dfrac{\Log(1/w)}{2w} \lp 1 + O\lp \dfrac{1}{\Log(1/w)} \rp \rp, \\
        F_2\lp G(1;q), G^{(1)}(1;q) \rp &= G^{(1)}(1;q) + G^2(1;q) = \dfrac{\Log^2(1/w)}{4w^2} \lp 1 + O\lp \dfrac{1}{\Log(1/w)} \rp \rp, \\
        F_3\lp G(1;q), G^{(1)}(1;q), G^{(2)}(1;q) \rp &= G^{(2)}(1;q)+3G^{(1)}(1;q)+G^{(1)}(1;q)G(1;q)+G^{3}(1;q) \\ &= \dfrac{\Log^3(1/w)}{8w^3} \lp 1 + O\lp \dfrac{1}{\Log(1/w)} \rp \rp,
    \end{align*}
    and by induction of $m$ we likewise obtain
    \begin{align*}
        F_m\lp G(1;q), \dots, G^{(m-1)}(1;q) \rp = \dfrac{\Log^m(1/w)}{2^m w^m} \lp 1 + O\lp \dfrac{1}{\Log(1/w)} \rp \rp.
    \end{align*}
    By Lemma \ref{L: No z}, and with the help of \cite[Lemma 2.3]{BringmannCraigMalesOno} we also have
    \begin{align*}
        \Log\lp q;q \rp_\infty^{-1} = \dfrac{\pi^2}{6w} + \dfrac{\Log(w)}{2} - \dfrac{\Log(2\pi)}{2} + O(w),
    \end{align*}
    from which we derive
    \begin{align*}
        P(1;q) = \lp q;q \rp_\infty^{-1} = \sqrt{\dfrac{w}{2\pi}} \exp\lp \dfrac{\pi^2}{6w} + O(w) \rp.
    \end{align*}
    We therefore conclude that
    \begin{align*}
        \mathcal A_m(q) &= F_m\lp G(z;q), G^{(1)}(z;q), \dots, G^{(m-1)}(z;q) \rp \cdot P(z;q) \\ &= \dfrac{\Log^m(1/w)}{2^m \sqrt{2\pi} w^{m - \frac 12}} e^{\frac{\pi^2}{6w}} \lp 1 + O\lp \dfrac{1}{\Log(1/w)} \rp \rp
    \end{align*}
    as desired.
\end{proof}

We now move on to identify the ``minor arc" behavior of $\mathcal A_m(q)$, which essentially says that $\mathcal A_m(q)$ is extremely small away from $q = 1$.

\begin{lemma} \label{L: Minor arc condition}
    As $|w| \to 0$ in the bounded cone $\frac{\pi}{2} - \delta \leq |\arg(w)| < \frac{\pi}{2}$, we have
    \begin{align*}
        \left| \mathcal A_m\lp q \rp \right| \ll_\delta \dfrac{|\Log(1/w)|^m}{|w|^{m-\frac 12}} e^{\lp \frac{\pi^2}{6} - \varepsilon \rp \frac 1x}
    \end{align*}
    for some constant $\varepsilon > 0$.
\end{lemma}

\begin{remark}
    This result is also easily demonstrated using the modularity of Dedekind's eta function, and can even be made explicit by combining this fact with Euler's pentagonal number theorem, see e.g. \cite{JacksonOtgonbayar}. We provide a different proof, which does not rely on modularity, since not all $q$-Pochhammer products are modular.
\end{remark}

\begin{proof}
    This follows quite quickly from Lemma \ref{L: No z}. In particular, the function $\mathcal A_m(q)$ is $\lp q;q \rp_\infty^{-1}$ multiplied by a function of only rational growth (as established in the proof of Proposition \ref{P: Major arc condition}). Since Lemma \ref{L: No z} establishes that $\lp q;q \rp_\infty^{-1}$ has the requisite exponential savings in open sets near roots of unity. By choosing a finite subcover of the unit circle, which can be done by compactness, such exponential savings can be extended to the whole unit circle by taking suprema.
    
    For a more detailed variation of this style of argument, see for instance the arguments such as those in Sections 4 and 5 in \cite{BridgesFrankeGarnowski}, which explains a more precise breakdown of exactly which open sets to choose (i.e. restricting to $k \leq n^\delta$ for suitable $\delta$ gives appropriate error terms). We do not replicate the exact arguments here since they are largely {\it mutatis mutandis}.
\end{proof}

\section{Proof of Theorem \ref{T: Main Theorem}} \label{S: Proof}

We now use the results of previous sections along with Wright's Circle Method to prove an asymptotic formula for $A_m(n)$, and thus also for $\mathbb{E}_m(n)$. We follow in particular the work Ngo and Rhoades \cite{NgoRhoades} and of Kim and Kim \cite{KimKim} in using the Circle Method to study partition statistics, though our asymptotic formulas differ from theirs.

By Cauchy's Theorem, we know that
\begin{align*}
    A_m(n) = \dfrac{1}{2\pi i} \int_{\mathcal C} \dfrac{\mathcal A_m(q)}{q^{n+1}} dq
\end{align*}
for any circle $\mathcal C$ centered at the origin with radius less than 1. Now, we choose $q = e^{-w}$, $w = \frac{\pi}{\sqrt{6n}}\lp 1 + iy \rp$. Set $\Delta = \cot(\delta) > 0$ for $0 < \delta < \frac{\pi}{2}$, noting that then $\Delta>0$. Choosing $\mathcal C$ to have radius $|q| = e^{-\pi/\sqrt{6n}}$, we may then write
\begin{align*}
    A_m(n) = \dfrac{1}{2\pi i} \int_{\mathcal C} \dfrac{\mathcal A_m(q)}{q^{n+1}} dq = M_m(n) + E_m(n),
\end{align*}
where we have the major arc integral
\begin{align*}
    M_m(n) = \dfrac{1}{2\pi} \lp \dfrac{\pi}{\sqrt{6n}} \rp \int_{|y| < \Delta} \mathcal A_m(e^{-w}) e^{nw} dy,
\end{align*}
and the minor arc integral
\begin{align*}
    E_m(n) = \dfrac{1}{2\pi} \lp \dfrac{\pi}{\sqrt{6n}} \rp \int_{\Delta \leq |y| < \sqrt{6n}} \mathcal A_m(e^{-w}) e^{nw} dy.
\end{align*}
By Proposition \ref{P: Major arc condition}, we may decompose $M_m(n)$ as
\begin{align*}
    M^{(1)}_m(n) &= \dfrac{1}{2^{m+1} \pi \sqrt{2\pi}} \lp \dfrac{\pi}{\sqrt{6n}} \rp \int_{|y|<\Delta} w^{-m+\frac 12} \Log^m(1/w) \exp\lp \dfrac{\pi^2}{6w} + nw \rp dy, \\
    M^{(2)}_m(n) &= M_m(n) - M^{(1)}_m(n).
\end{align*}
Then we have
\begin{align*}
    M_m^{(1)}(n) &= \dfrac{1}{2^{m+1}\pi\sqrt{2\pi}} \lp \dfrac{\pi}{\sqrt{6n}} \rp^{\frac 32 - m} \int_{-\Delta}^\Delta \lp 1 + iy \rp^{-m+\frac 12} \lp \log\lp \dfrac{\sqrt{6n}}{\pi} \rp - \Log\lp 1 + iy \rp \rp^m \\ &\cdot \exp\lp \dfrac{\pi}{\sqrt{6}} \sqrt{n} \lp \dfrac{1}{1+iy} + 1 + iy \rp \rp dy,
\end{align*}
which by expanding via the binomial theorem obtains
\begin{align*}
    &M_m^{(1)}(n) = \dfrac{1}{2^{m+1}\pi\sqrt{2\pi}} \lp \dfrac{\pi}{\sqrt{6n}} \rp^{\frac 32 - m} \\ &\cdot\sum_{j=0}^m \binom{m}{j} \log^{m-j}\lp \dfrac{\sqrt{6n}}{\pi} \rp \int_{-\Delta}^\Delta \dfrac{(-1)^j \Log^j\lp 1 + iy \rp}{\lp 1 + iy \rp^{m-\frac 12}} \exp\lp \dfrac{\pi}{\sqrt{6}} \sqrt{n} \lp \dfrac{1}{1+iy} + 1 + iy \rp \rp dy.
\end{align*}

For further convenience, set $f(y) := \frac{1}{1+iy} + 1 + iy$. We note that $\mathrm{Re}\lp f(y) \rp \leq 2 = f(0)$ for $|y| \leq \Delta$, $f^\prime(0) = 0$, and $f^{\prime\prime}(0) = -2$. Thus, by the saddle point method, i.e. Theorem \ref{T: Saddle point method}, we have 
\begin{align*}
    \int_{-\Delta}^\Delta \lp 1 + iy \rp^{-m+\frac{1}{2}} \exp\lp \dfrac{\pi}{\sqrt{6}} \sqrt{n} \lp \dfrac{1}{1+iy} + 1 + iy \rp \rp dy = \lp \dfrac{6}{n} \rp^{\frac 14} e^{\pi\sqrt{\frac{2n}{3}}} \lp 1 + O\lp n^{-\frac 14} \rp \rp.
\end{align*}

Because $\left| \Log\lp 1 + iy \rp \right| \leq |y|$ and $\mathrm{Re}\lp f(y) \rp \leq 2 - \frac{y^2}{1+\Delta^2}$ for $|y| \leq \Delta$, it is straightforward to show that for $j \geq 1$ that
\begin{align*}
    \int_{-\Delta}^\Delta \dfrac{\Log^j\lp 1 + iy \rp}{\lp 1 + iy \rp^{m-\frac 12}} \exp &\lp \dfrac{\pi}{\sqrt{6}} \sqrt{n} \lp \dfrac{1}{1+iy} + 1 + iy \rp \rp dy \\ &\ll_\Delta e^{\pi\sqrt{\frac{2n}{3}}} \int_{-\Delta}^\Delta |y|^j \exp\lp -\dfrac{\pi \sqrt{n} y^2}{\sqrt{6}\lp 1 + \Delta^2 \rp} \rp dy \\ &\ll_\Delta e^{\pi\sqrt{\frac{2n}{3}}} n^{-\frac j2} \int_{-\infty}^\infty |y|^j \exp\lp -\dfrac{\pi y^2}{\sqrt{6}\lp 1 + \Delta^2 \rp} \rp dy \\ &\ll_\Delta e^{\pi\sqrt{\frac{2n}{3}}} n^{-\frac j2}.
\end{align*}
It is similarly straightforward to show that $M_m^{(2)}(n)$ grows at least $\frac{1}{\log(n)}$ slower than $M_m^{(1)}(n)$, since $M_m^{(1)}(n)$ characterizes the main asymptotic contribution of $\mathcal A_m(q)$ up to $O\lp \frac{1}{\Log(w)} \rp$ by Proposition \ref{P: Major arc condition}. Thus, we have
\begin{align*}
    M_m(n) &= \dfrac{1}{2^{m+1}\pi\sqrt{2\pi}} \lp \dfrac{\pi}{\sqrt{6n}} \rp^{\frac 32 - m} \lp \dfrac{6}{n} \rp^{\frac 14} \log^m\lp \dfrac{\sqrt{6n}}{\pi} \rp e^{\pi\sqrt{\frac{2n}{3}}} \lp 1 + O\lp \dfrac{1}{\log(n)} \rp \rp
    \\ &= \dfrac{n^{\frac{m-2}{2}}}{2^{m+2}\sqrt{3}} \lp \dfrac{\sqrt{6}}{\pi} \rp^m \log^m\lp \dfrac{\sqrt{6n}}{\pi} \rp e^{\pi\sqrt{\frac{2n}{3}}} \lp 1 + O\lp \dfrac{1}{\log(n)} \rp \rp.
\end{align*}
The result then follows from the minor arc bound, derived from Lemma \ref{L: Minor arc condition}, that
\begin{align*}
    \left| E_m(n) \right| &\ll_\delta \dfrac{1}{\sqrt{n}} \int_{\Delta \leq |y| < \sqrt{6n}} |w|^{-m+\frac 12} \left| \Log(w) \right|^{m+1} \exp\lp n x + \lp \dfrac{\pi^2}{6} - \varepsilon \rp \dfrac{1}{x} \rp dy \\ &\ll \log^{m+1}(n) \exp\lp \pi \sqrt{\dfrac{2n}{3}} - \dfrac{\sqrt{6} \varepsilon}{\pi} \sqrt{n} \rp,
\end{align*}
which now proves that 
\begin{align*}
    A_m(n) = \dfrac{n^{\frac{m-2}{2}}}{2^{m+2}\sqrt{3}} \lp \dfrac{\sqrt{6}}{\pi} \rp^m \log^m\lp \dfrac{\sqrt{6n}}{\pi} \rp e^{\pi\sqrt{\frac{2n}{3}}} \lp 1 + O\lp \dfrac{1}{\log(n)} \rp \rp.
\end{align*}
after we recall the Hardy--Ramanujan asymptotic formula for partitions \cite{HardyRamanujan}, namely
\begin{align*}
    p(n) \sim \dfrac{1}{4n\sqrt{3}} e^{\pi \sqrt{\frac{2n}{3}}}.
\end{align*}

\section{Further Discussions} \label{S: Further Discussion}

\subsection{MacMahon's partition analysis}

MacMahon's partition analysis provides a powerful tool to compute the multi-variable generating function for all kinds of partitions, especially when the restriction between different parts can be expressed by linear inequalities. In combination with the Circle Method and the other analytic techniques used already in this paper, we suggest that these methods in tandem should be of substantial use in computing asymptotics for joint moments in probability distributions, and therefore for computing full joint distributions of partition statistics, perhaps even many at the same time.

In Section $2.1$ we have seen the product form for the set of all partitions. However, in practice, it is more custom to get the series form by considering partitions with fixed length. Taking the set of all partitions as an example, we would have
\begin{align*}
\sum_{\substack{\lambda\in\mathcal{P}\\\ell(\lambda)=n}}x^{\lambda}=&\sum_{a_1\geq a_2\geq\cdots\geq a_n>0}x_1^{a_1}x_2^{a_2}\cdots x_n^{a_n}\\
    =&\underset{\geq}{\Omega}\sum_{a_1,a_2,\ldots,a_n\geq0}x_1^{a_1}x_2^{a_2}\cdots x_n^{a_n}\lambda_1^{a_1-a_2}\lambda_2^{a_2-a_3}\cdots\lambda_{n-1}^{a_{n-1}-a_n}\lambda_n^{a_n-1}\\
    =&\underset{\geq}{\Omega}\frac{\lambda_n^{-1}}{(1-x_1\lambda_1)(1-\frac{x_2\lambda_2}{\lambda_1})\cdots(1-\frac{x_{n-1}\lambda_{n-1}}{\lambda_{n-2}})(1-\frac{x_n\lambda_n}{\lambda_{n-1}})}.
\end{align*}
By eliminating the $\lambda$'s in this expression, we have
$$\sum_{\substack{\lambda\in\mathcal{P}\\\ell(\lambda)=n}}x^{\lambda}=\frac{X_n}{(1-X_1)(1-X_2)\cdots(1-X_n)}.$$
Now by adding up all the nonnegative integers $n$ and compare with what we get in the proof of Lemma \ref{generating function}, we have
$$\sum_{\lambda\in\mathcal{P}}x^{\lambda}=\sum_{n=0}^{\infty}\frac{X_n}{\prod_{i=1}^{n}(1-X_i)}=\prod_{n=1}^{\infty}\frac{1}{1-X_n}.$$
We have already seen how this help us tracing the alternating sum of a partition, meanwhile, we can consider other statistics. For instance, let $L(\lambda)$ be the largest part of $\lambda$, and note that the largest part is always indicated by $x_1$, with proper substitutions we have
\begin{align*}
\sum_{\lambda\in\mathcal{P}}x^{L(\lambda)}z^{a(\lambda)}q^{|\lambda|}=&1+\sum_{n=1}^{\infty}\frac{xq^{2n}}{(xzq;q^2)_n(xq^2;q^2)_n}+\sum_{n=0}^{\infty}\frac{xzq^{2n+1}}{(xzq;q^2)_{n+1}(xq^2;q^2)_n}\\
=&\frac{1}{(xzq;q^2)_{\infty}(xq^2;q^2)_{\infty}}.    
\end{align*}
Thus we get the three-variable generating function and the sum-product identity automatically. Moreover, since the series form is based on the counting of the length, we can further get
$$\sum_{\lambda\in\mathcal{P}}x^{L(\lambda)}y^{\ell(\lambda)}z^{a(\lambda)}q^{|\lambda|}=1+\sum_{n=1}^{\infty}\frac{xy^{2n}q^{2n}}{(xzq;q^2)_n(xq^2;q^2)_n}+\sum_{n=0}^{\infty}\frac{xzy^{2n+1}q^{2n+1}}{(xzq;q^2)_{n+1}(xq^2;q^2)_n},$$
which is a four-parameter generating function for the set of all partitions. Unfortunately, the length of the partition is not directly indicated in the product form, but the series form is always guaranteed by the partition analysis. In addition, with the help of the theory of hypergeometric series, one might still convert the series into a product. The same method can be applied to the set of strict partitions, that is, partitions without repeating parts. Let $\mathcal{S}$ be the set of all such partitions, then following the same steps, we have
\begin{align*}
\sum_{\substack{\lambda\in\mathcal{S}\\\ell(\lambda)=n}}x^{\lambda}=&\sum_{a_1>a_2>\cdots>a_n>0}x_1^{a_1}x_2^{a_2}\cdots x_{n}^{a_n}\\
=&\sum_{a_1,a_2,\ldots,a_n\geq0}x_1^{a_1}x_2^{a_2}\cdots x_{n}^{a_n}\lambda_1^{a_1-a_2-1}\lambda_2^{a_2-a_3-1}\cdots\lambda_{n-1}^{a_{n-1}-a_n-1}\lambda_{n}^{a_n-1}\\
=&\frac{\lambda_1^{-1}\lambda_2^{-1}\cdots\lambda_{n}^{-1}}{(1-x_1\lambda_1)(1-\frac{x_2\lambda_2}{\lambda_1})\cdots(1-\frac{x_{n-1}\lambda_{n-1}}{\lambda_{n-2}})(1-\frac{x_n\lambda_n}{\lambda_{n-1}})},
\end{align*}
which is the crude form for the generating function. By eliminating the $\lambda$'s and adding up all the nonnegative $n$, we have
$$\sum_{\lambda\in\mathcal{S}}x^{\lambda}=\sum_{n=0}^{\infty}\frac{x_1^{n}x_2^{n-1}\cdots x_{n}}{(1-x_1)(1-x_1x_2)\cdots(1-x_1x_2\cdots x_n)}=\sum_{n=0}^{\infty}\prod_{i=1}^{n}\frac{Xi}{1-X_i}.$$
Now with the same substitution, we have
\begin{align*}
\sum_{\lambda\in\mathcal{S}}x^{L(\lambda)}y^{\ell(\lambda)}z^{a(\lambda)}q^{|\lambda|}&=\sum_{n=0}^{\infty}\frac{x^{2n}y^{2n}z^{n}q^{\binom{2n+1}{2}}}{(xzq;q^2)_n(xq^2;q^2)_n}+\sum_{n=0}^{\infty}\frac{x^{2n+1}y^{2n+1}z^{n+1}q^{\binom{2n+2}{2}}}{(xzq;q^2)_{n+1}(xq^2;q^2)_n}.
\end{align*}
Here we have treated the generating function for ordinary partitions and strict partitions. More applications of MacMahon's partition analysis can be found in the references we mentioned in the introduction. Meanwhile, the method can also be applied to the study of $2$-dimensional partitions, please refer to \cite{Andrews1,AndrewsPaule1} for more details.
\subsection{Implications for distributions}

To close out the paper, we outline the basic features of our approach in this work to studying distributions for partition statistics and its particular merits. Methods similar to this have been implemented by other authors as well; we note in particular the implementation by Bringmann, Mahlburg and Rhoades \cite{BringmannMahlburgRhoades}. In their work, they study the {\it rank} and {\it crank} statistics for partitions, which we call $r(\lambda)$ and $c(\lambda)$, respectively\footnote{Since we do not perform any technical calculations with these, we do not require their definitions.}. We let $M(m,n)$ and $N(m,n)$ be the number of partitions of $n$ with crank $m$ or rank $m$, respectively. The generating functions for these statistics are given by
\begin{align*}
    C(x;q) &:= \sum_{n \geq 0} M(m,n) x^m q^b = \prod_{n \geq 1} \dfrac{1 - q^n}{\lp 1 - xq^n \rp \lp 1 - x^{-1} q^n \rp} = \dfrac{1-x}{\lp q;q \rp_\infty} \sum_{n \geq 0} \dfrac{(-1)^n q^{\frac{n(n+1)}{2}}}{1-xq^n}, \\ 
    R(x;q) &:= \sum_{n \geq 0} N(m,n) x^m q^n = \sum_{n \geq 0} \dfrac{q^{n^2}}{\lp xq;q \rp_n \lp x^{-1}q;q \rp_n} = \dfrac{1-x}{\lp q;q \rp_\infty} \sum_{n \in \Z} \dfrac{(-1)^n q^{\frac{n(3n+1)}{2}}}{1-xq^n}.
\end{align*}
Their work is devoted to the computation of asymptotics for the moments
\begin{align*}
    M_k(n) &:= \sum_{\lambda \vdash n} c(\lambda)^k = \sum_{m \in \Z} m^k M(m,n), \\
    N_k(n) &:= \sum_{\lambda \vdash n} r(\lambda)^k = \sum_{m \in \Z} m^k N(m,n).
\end{align*}
Based on this work, Diaconis, Janson and Rhoades \cite{DiaconisJansonRhoades} identify the distribution laws for these statistics. This result follows by identifying these asymptotics as giving rise to particular moment generating functions. Philosophically, the takeaway is that {\it method of moments} (i.e. the computation of moments of a partition statistic) yields the full asymptotic distribution of a statistic. This is a well-known principle in statistics (see e.g. \cite{Curtiss}), and our approach fits exactly this model and identifies the distribution of $a(\lambda)$ directly by Theorem \ref{T: Main Theorem}.

We hope very strongly that this method will be fruitful in the theory of partitions, possibly moreso than other existing methods, with the following evidence. Consider a statistic $s(\lambda)$ on partitions, with the two-variable generating function
\begin{align*}
    S(z;q) := \sum_{\lambda \in \mathcal P} z^{s(\lambda)} q^{|\lambda|}.
\end{align*}
Then it is easy to see that, for the differential operator $\partial_z := z \frac{\partial}{\partial z}$, that the $k$th moment generating function
\begin{align*}
    S_k(z;q) = \sum_{\lambda} s(\lambda)^k q^{|\lambda|}
\end{align*}
is given by
\begin{align*}
    S_k(z;q) = \partial_z^k S(z;q) |_{z=1}.
\end{align*}
We thus obtain all moments by taking $z$-derivatives. It is probable that such operations will be useful in the theory of partitions. As seen in \cite{BringmannMahlburgRhoades}, any case where the generating function possesses modular or Jacobi structure will permit the computation of precise asymptotics for moments in all cases. Even when this is not true, there are many techniques (such as Euler--Maclaurin summation as used in our work) which can also be used to calculate asymptotics for these derivatives. For example, looking at the case $\lp zq^a;q^b \rp_\infty^{-1}$ in this paper, we see
\begin{align*}
    \partial_z \lp zq^a; q^b \rp_\infty^{-1} = \lp \sum_{k \geq 0} \dfrac{z q^{a+bk}}{1 - zq^{a+bk}} \rp \lp zq^a;q^b \rp_\infty^{-1}.
\end{align*}
Thus, there is a natural recursive structure to large families of partition generating functions which leads to the need for computation of asymptotics for {\it Lambert series}, which take roughly the form
\begin{align*}
    \sum_{n \geq 0} \dfrac{a_n P_k\lp q^n \rp}{\lp 1 - q^n \rp^k}
\end{align*}
for $k \geq 1$ and $P_k(q)$ a polynomial. In very many situations, asymptotics for such sums can be computed using Euler--Maclaurin summation or modular/Jacobi techniques. Asymptotics of this sort for generating functions would then lead to asymptotic formulas for the moments of $s$ using Wright's Circle Method \cite{NgoRhoades} or Ingham's Tauberian Theorem \cite{BringmannJenningsShafferMahlburg}. We therefore propose the following proposed direction of study for partition theory.

\begin{question}
    We pose the following questions:
    \begin{itemize}
        \item For broad families of partition statistics $s : \mathcal P \to \C$ with generating function $S(z;q)$, consider the asymptotic structure of the derivatives $\partial_z S(z;q)$. What is the recursive and asymptotic structure of such derivatives?
        \item What does this recursive structure dictate about distribution laws for partition statistics?
    \end{itemize}
\end{question}

We also hope that techniques from MacMahon's partition analysis can help produce an easy-to-access theory of joint distributions.

\begin{question}
    Can one prove, from multi-variate generating functions coming from MacMahon's partition analysis, joint distributions laws for pairs of interesting partition statistics?
\end{question}

We hope these questions will motivate further study between the interrelationship between common methods for constructing generating functions in partition theory and corresponding asymptotic methods. Based on this work as well as several other works computing asymptotics for various statistics which have been cited already, we propose a potential framework which might describe the main structures that have arisen.

\begin{question}
    How can one characterize families of partition statistics that necessarily obtain the same basic type of distribution law?
\end{question}

We also suggest that this method could be very fruitful for computing distributions of statistics arising from partition functions in statistical mechanics (see for instance \cite{}). We leave such observations or implementations for future work.

\end{document}